\documentclass{amsproc}

\usepackage[T1]{fontenc}
\usepackage{fancyhdr}

\usepackage{amsfonts}
\usepackage{amssymb,amsthm,amsmath,amstext,amsxtra}
\usepackage{relsize}
\usepackage{bm}       % allows bold italic letters in math mode
\usepackage{bbm}       % Blackboard variants of Computer Modern fonts
\usepackage{mathtools} % allows more extendible arrows
\usepackage{colonequals}
\usepackage{numprint} % format long numbers
\usepackage[shortlabels]{enumitem}
\usepackage{chngcntr}

\usepackage{tikz}
\usepackage{tikz-cd}
\usepackage{xcolor}
\usepackage{color}
\usepackage[pdfencoding=auto, psdextra, breaklinks=true]{hyperref}
\definecolor{mylinkcolor}{rgb}{0.5,0.0,0.0}
\definecolor{myurlcolor}{rgb}{0.0,0.0,0.7}
\hypersetup{
 colorlinks,
 urlcolor=myurlcolor,
 citecolor=myurlcolor,
 linkcolor=mylinkcolor,
 breaklinks=true,
 pdfauthor={Shiva Chidambaram}
}
\usepackage[capitalize]{cleveref}

\usepackage{csquotes} % removes one warning of biblatex
\usepackage[style=alphabetic, backend=biber, sorting=nyt, doi=false, isbn=false,url=false, hyperref,backref,backrefstyle=none, maxnames=20, maxalphanames=20]{biblatex}
\newbibmacro{string+doiurlisbn}[1]{%
  \iffieldundef{doi}{%
    \iffieldundef{url}{%
      \iffieldundef{isbn}{%
        \iffieldundef{issn}{%
          #1%
        }{%
          \href{https://books.google.com/books?vid=ISSN\thefield{issn}}{#1}%
        }%
      }{%
        \href{https://books.google.com/books?vid=ISBN\thefield{isbn}}{#1}%
      }%
    }{%
      \href{\thefield{url}}{#1}%
    }%
  }{%
    \href{https://dx.doi.org/\thefield{doi}}{#1}%
  }%
}
\DeclareFieldFormat{title}{\usebibmacro{string+doiurlisbn}{\mkbibemph{#1}}}
\DeclareFieldFormat[article,incollection,inproceedings,unpublished,misc,book]{title}%
    {\usebibmacro{string+doiurlisbn}{\mkbibquote{#1}}}
\DefineBibliographyStrings{english}{%
  backrefpage = {$\uparrow$},%
  backrefpages = {$\uparrow$}%
}
\usepackage{xurl}
\usepackage{breakurl}

\counterwithin*{equation}{section}

\def\Gal{\mathrm{Gal}}

\def\SL{\mathrm{SL}}

\DeclareMathOperator{\frob}{Frob}
\DeclareMathOperator{\proj}{pr}
\DeclareMathOperator{\tr}{Tr}
\DeclareMathOperator{\Hom}{Hom}

\DeclareMathOperator{\Norm}{Nm}
\DeclareMathOperator{\Ram}{Ram}
\DeclareMathOperator{\coker}{coker}
\DeclareMathOperator{\res}{res}
\DeclareMathOperator{\GL}{GL}
\DeclareMathOperator{\GSp}{GSp}
\DeclareMathOperator{\Sp}{Sp}

\DeclareMathOperator{\Aut}{Aut}

\newcommand{\MA}{\mathcal{A}}

\newcommand{\F}{\mathbb{F}}
\newcommand{\N}{\mathbb{N}}
\newcommand{\Q}{\mathbb{Q}}

\newcommand{\Z}{\mathbb{Z}}
\newcommand{\calX}{\mathcal{X}}

\newcommand{\Qbar}{\overline{\Q}}
\newcommand{\vv}{\mathbf{v}}

\newcommand\rhobar{\overline{\rho}}

\DeclareFontFamily{U}{wncy}{}
\DeclareFontShape{U}{wncy}{m}{n}{<->wncyr10}{}
\DeclareSymbolFont{mcy}{U}{wncy}{m}{n}
\DeclareMathSymbol{\Sh}{\mathord}{mcy}{"58}

\DeclarePairedDelimiter\floor{\lfloor}{\rfloor}

\numberwithin{equation}{section}
\newtheorem{prop}[equation]{Proposition}
\Crefname{prop}{Proposition}{Propositions}
\newtheorem{thm}[equation]{Theorem}
\Crefname{thm}{Theorem}{Theorems}
\newtheorem{lem}[equation]{Lemma}
\Crefname{lem}{Lemma}{Lemmas}

\Crefname{cor}{Corollary}{Corollaries}
\newtheorem{question}[equation]{Question}
\Crefname{question}{Question}{Questions}

\theoremstyle{definition}
\newtheorem{defn}[equation]{Definition}

\title[Mod-$p$ Galois representations not arising from abelian varieties]{Mod-$p$ Galois representations not arising from abelian varieties}

\author[S. Chidambaram]{Shiva Chidambaram}
\address{Department of Mathematics, Massachusetts Institute of Technology, Cambridge, MA 02139-4307, USA}
\email{shivac@mit.edu}
\urladdr{\url{https://math.mit.edu/~shivac}}


\begin{document}

\begin{abstract}
It is known that any Galois representation $\rho : G_{\Q} \rightarrow \GL(2,\F_p)$ with determinant equal to the cyclotomic character, arises from the $p$-torsion of an elliptic curve over $\Q$, if and only if $p \leq 5$. In dimension $g = 2$, when $p \le 3$, it is again known that any Galois representation valued in $\GSp(4,\F_p)$ with cyclotomic similitude character arises from an abelian surface. In this paper, we study this question for all primes $p$ and dimensions $g \ge 2$. When $g \ge 2$ and $(g,p) \neq (2,2)$, $(2,3)$, $(3,2)$, we prove the existence of a Galois representation over $\Q$ valued in~$\GSp(2g,\F_p)$ with cyclotomic similitude character, that cannot arise as the~$p$-torsion representation of any~$g$-dimensional abelian variety over~$\Q$.
\end{abstract}

\maketitle

\renewcommand{\thefootnote}{\fnsymbol{footnote}} 
\renewcommand{\thefootnote}{\arabic{footnote}}
\section{Introduction}
\label{sec:intro}

Let $g \geq 1$ and $p$ be a prime. 
Let $\MA_g(p)$ be the Siegel modular variety which is the moduli space of principally polarized abelian varieties of dimension~$g$ with full level $p$ structure. When $g = 1$, this is exactly the modular curve $X(p)$ parametrizing elliptic curves with full level $p$ structure. Using Riemann-Hurwitz formula, one knows that the curve $X(p)$ has genus $0$ if and only if $p \leq 5$. In general, the variety~$\MA_g(p)$ is geometrically rational if and only if $(g,p) = (1,2),$ $(1,3),$ $(1,5),$ $(2,2),$ $(2,3),$ $(3,2)$ \cite{Hulek}.

In all three genus $1$ cases above, the modular curve $X(p) = \MA_1(p)$ is in fact rational over $\Q$, and furthermore all of its twists $\MA_1(\rho)$ corresponding to two dimensional mod-$p$ Galois representations $\rho$ with cyclotomic determinant, are also rational over $\Q$ \cite{RS}. That is, if $p = 2, 3$ or $5$, and $\rho : G_{\Q} \rightarrow \GL(2,\F_p)$ is such that $\det \rho$ is the mod-$p$ cyclotomic character~$\chi_p$, then the Galois representation $\rho$ arises from an elliptic curve over $\Q$, and in fact from infinitely many elliptic curves. 
When $g = 2$ and $p = 2$ or $3$, the moduli spaces~$\MA_2(p)$ and
their twists are all unirational, as explained in~\cite[Lemma 10.2.4]{BCGP} and~\cite[\S4.4]{CCR}, even if they are not always rational over $\Q$ \cite{CC}. Hence, in these cases, all representations $\rho: G_{\Q} \rightarrow \GSp(4,\F_p)$ with cyclotomic similitude character do arise from abelian surfaces over $\Q$.

In this paper, we consider the cases where $g \geq 2$ and $\MA_g(p)$ is not geometrically rational.
The main theorem we prove is:
\begin{thm}
\label{thm:main-theorem}
    Let $g \geq 2$ and $p$ be a prime number. Suppose $(g,p) \neq (2,2)$, $(2,3)$ or $(3,2)$ so that $\MA_g(p)$ is not geometrically rational. Then there exists a Galois representation $\rho: G_{\Q} \rightarrow \GSp(2g,\F_p)$ with cyclotomic similitude character such that $\rho$ does not arise as the $p$-torsion representation of any 
abelian variety over $\Q$.
\end{thm}

This leaves the question open only in the case $(g,p) = (3,2)$.
\begin{question}
    Given any Galois representation~$\rho : G_{\Q} \rightarrow \Sp(6,\F_2)$, is it true that $\rho$ arises as the $2$-torsion representation $\rhobar_{A,2}$ of some abelian threefold $A$ over $\Q$?
\end{question}

The analogous statement for primes $p > 5$ in dimension $g = 1$ was proved independently in \cite{Cal06} and \cite{Die04}, showing the existence of non-elliptic mod $p$ Galois representations. For $p \geq 11$, \cite{Cal06} constructs the desired representations by considering certain weight $2$ Hecke eigenforms with coefficient field not equal to $\Q$. For $p = 7$, this method does not work, and the desired non-elliptic representation is constructed explicitly, with image contained in the normalizer of the non-split Cartan subgroup of $\GL(2,\F_7)$. Our proof of \cref{thm:main-theorem} involves generalizing this explicit construction to higher genus situations.

Let $\rho : G_{\Q} \rightarrow \GSp(2g,\F_p)$ be a mod-$p$ Galois representation with cyclotomic similitude. For any prime $\ell \neq p$, let $I_{\ell} = \Gal(\Qbar_{\ell}|\Q_{\ell}^{ur})$ denote the absolute inertia group at $\ell$.  image under $\rho$ of the inertia group at $\ell$. Suppose that $\rho$ arises from an abelian variety $A$ over $\Q$, i.e., $\rho \simeq \rhobar_{A,p}$. Then Raynaud’s inertial criteria for semistable reduction \cite[\S 7.4, Thm 6]{BLR90} implies that there exists a constant $K_g$ only depending on $g$ such that the prime-to-$p$ part of $\#\rho(I_{\ell})$ 
divides $K_g$. On the other hand, we can obtain using Zsigmondy’s theorem \cite{Zsigmondy}, that there exists an integer $q>1$ coprime to $p$ such that~$q | \#\GSp(2g,\F_p)$ and $q \nmid K_g$. If we could construct a Galois representation $\rho$ such that $\rho(I_{\ell})$ had order~$q$, we could immediately deduce that $\rho$ does not come from any $g$-dimensional abelian variety. Thus this becomes an instance of
the Inverse Galois Problem with local conditions.

Here is an outline of the paper. In \cref{sec:abvar}, we recall properties related to semistable reduction of abelian varieties, and prove the existence of the constant $K_g$ controlling the prime-to-$p$ part of order of inertia at primes $\ell \ne p$. We show that the only primes dividing $K_g$ are the ones less than or equal to $2g+1$. 
In \cref{sec:cartan}, we describe a certain metabelian subgroup $N$ inside the normalizer of a non-split Cartan subgroup of~$\GSp(2d,\F_p)$. It is an extension of the cyclic group $\Gal(\F_{p^{2d}}|\F_p) \simeq \Z/2d$ by a cyclic group of order $(p^d+1)(p-1)$ (a non-split Cartan coming from the field extension $\F_{p^{2d}}|\F_{p^d}$.) The restriction of the similitude character to $N$ remains surjective. In \cref{sec:embed}, we study the embedding problem
\begin{equation}
\begin{tikzcd}
    & & & G_{\Q} \arrow{d}{\phi} \arrow[dashed]{dl}[swap]{?}& \\
    0 \arrow{r} & \left[ N,N \right] \arrow{r} & N \arrow{r} & N^{ab} \arrow{r} & 0
\end{tikzcd}
\end{equation}
for homomorphisms $\phi$, such that the composition of $\phi$ with the similitude character $N^{ab} \rightarrow \F_{p}^{\times}$ is equal to the mod-$p$ cyclotomic character $\chi_p$. We choose $\phi$ carefully so that all ramification is tame, and then proceed to show in Sections \ref{subsec:infty}-\ref{subsec:N2} that all local obstructions to the embedding problem vanish. Using Galois cohomological machinery related to the Grunwald-Wang theorem \cite[Chapter IX]{Neukirch2008}, we show that global obstructions vanish as well, and hence that $\phi$ can be lifted to a proper solution $\tilde{\phi} : G_{\Q} \twoheadrightarrow N$.
In \cref{sec:proof-main-theorem}, we finish the proof of \cref{thm:main-theorem}.
We show that by twisting $\tilde{\phi}$ suitably, we can obtain representations $\rho$ 
such that $\rho(I_{\ell}) \subset [N,N]$ has prime power order $q$ not dividing $K_g$. By allowing ourselves to consider reducible representations landing inside $\GSp(2d,\F_p) \subset \GSp(2g,\F_p)$ for $d \leq g$, we can deal with all cases except $(g,p) = (3,3)$ using this approach. We deal with the exceptional case explicitly, 
by producing a representation whose image in $\GSp(6,\F_3)$ has order $78$, with $\rho(I_{\ell})$ being the unique cyclic subgroup of order $13$.

\section{Semistable reduction of abelian varieties}
\label{sec:abvar}

Let $X$ be an abelian variety of dimension $g$ defined over a field $F$. Let $v$ be a discrete valuation on $F$, and $\calX$ denote the Neron model of $X$ at $v$. Let $\ell$ be the residue characteristic of $v$.

\begin{defn}
    $X$ is said to have good reduction at $v$, if the identity component of the special fiber of $\calX$ is an abelian variety.
\end{defn}
\begin{defn}
    $X$ is said to have semistable reduction at $v$, if the identity component of the special fiber of $\calX$ is an extension of an abelian variety by an affine torus.
\end{defn}

Let $I_v$ denote the absolute inertia group at the finite prime $v$ of $F$. For a rational prime $p$, let $X[p]$ and $T_p(X)$ denote the $p$-torsion subgroup and the $p$-adic Tate module of $X$ respectively. Then we have the following critera for semistable reduction in terms of inertial action on $T_p(X)$ and $X[p]$ \cite[Prop 3.5, 4.7]{Gro72}, \cite[\S 7.4, Thm 6]{BLR90}.

\begin{thm}[Grothendieck]
	\label{inertcrit1}
	Let $p \neq \ell$ be a prime. Then the following are equivalent.
	\begin{enumerate}
		\item $X$ has semistable reduction at $v$.
		\item $I_v$ acts unipotently on the Tate module $T_p(X)$.
	\end{enumerate}
\end{thm}

\begin{thm}[Raynaud]
	\label{inertcrit2}
	Let $m \geq 3$ be an integer not divisible by $\ell$, and suppose that all the $m$-torsion points of $X$ are defined over an extension of $F$ unramified at $v$. Then $X$ has semistable reduction at $v$.
\end{thm}

For $g \geq 2$, let $K_g$ denote the gcd over all odd primes $r$ of the cardinality of $\GSp(2g,\F_r)$. Before we get to the inertial condition described in \cref{sec:intro}, we prove a few lemmas about $K_g$, wherein we make repeated use of Dirichlet's theorem about the infinitude of primes in arithmetic progression. Let $\nu_p$ denote the $p$-adic valuation function normalized so that $\nu_p(p) = 1$.

\begin{lem}
\label{Kglemma1}
All primes dividing $K_g$ are less than or equal to $2g+1$. Further, if $g \geq 2$ and $p$ is a prime such that $2 < p \leq 2g+1$, then $\nu_p (K_g) < g^2$.
\begin{proof}
    Let us first recall the formula
    \begin{align}
    \label{eqn:orderofGSp}
        \# \GSp(2g,\F_r) = (r-1) r^{g^2} \prod_{i=1}^g (r^{2i}-1).
    \end{align}
    Let $p > 2g+1$ be a prime. Choose a primitive root $a \in \F_p^{\times}$ and let $r \equiv a \pmod p$ be a prime. Then $r^{2i} \not\equiv 1 \pmod p$ for $1 \leq i \leq g$. This shows that $p$ does not $\# \GSp(2g,\F_r)$, and hence $p$ does not divide $K_g$.
    
    For the second part, choose a primitive root $a \in (\Z/p^{g^2})^{\times}$ and let $r \equiv a \pmod {p^{g^2}}$ be a prime. Then, for $n \leq g^2$, the order of $r \in (\Z/p^n)^{\times}$ is $p^{n-1}(p-1)$. So $p^n$ divides a term $r^{2i}-1$ in \cref{eqn:orderofGSp} if and only if $p^{n-1}(p-1)$ divides $2i$. Using this observation we count the powers of $p$ to get that
    \begin{align*}
        \nu_p(\# \GSp(2g,\F_r)) & = \nu_p\left((r^2-1)(r^4-1) \cdots (r^{2g}-1)\right) \\
        & = \floor*{\frac{2g}{p-1}} + \floor*{\frac{2g}{p(p-1)}} + \floor*{\frac{2g}{p^2(p-1)}} + \cdots \\
        & \leq \floor*{\frac{2g}{2}} + \floor*{\frac{2g}{4}} + \floor*{\frac{2g}{8}} + \cdots < 2g \leq g^2
    \end{align*}
    since $g \geq 2$. Therefore, $\nu_p(K_g) < g^2$.
\end{proof}
\end{lem}

\begin{lem}
\label{Kglemma2}
For any $M>2$, the gcd of $\# \GSp(2g,\F_r)$ as $r$ ranges over primes bigger than $M$ is equal to $K_g$.
\begin{proof}
    Let $L$ denote the gcd of $\# \GSp(2g,\F_r)$ as $r$ ranges over primes bigger than $M$. Following the argument in the proof of \cref{Kglemma1}, no prime greater than $2g+1$ divides $L$.
    Let $p \leq 2g+1$ be a prime and suppose $\nu_p(K_g) = n$. Then, there exists some prime $r \neq 2$ such that $\nu_p(\# \GSp(2g,\F_r)) = n$. If $r > M$, it is immediate that $\nu_p(L) = n$ as well. Suppose $r \leq M$. By the second part of \cref{Kglemma1}, we know that $n < g^2$ and hence $r \neq p$ since $\nu_p(\# \GSp(2g,\F_p)) = g^2$. Choose a prime $\ell > M$ such that $\ell \equiv r \pmod{p^{n+1}}$. Then $\# \GSp(2g,\F_l) \equiv \# \GSp(2g,\F_r) \pmod{p^{n+1}}$, showing that $\nu_p(\# \GSp(2g,\F_l)) = n$. This shows that $\nu_p(L) = n$ as well. Hence, $K_g = L$ which is what we want.
\end{proof}
\end{lem}

\begin{prop}
	\label{inertialorderdividesKg}
	Let $p \neq \ell$ be a prime and let $\rho : G_F \rightarrow \Aut(X[p])$ denote the $p$-torsion representation coming from the $g$-dimensional abelian variety $X$. Then, the prime-to-$p$ part of $\# \rho(I_v)$ divides $K_g$.
	\begin{proof}
		Suppose $X$ admits a polarization $X \rightarrow X^{\vee}$ of degree $M$. Thus for primes $r>M$, the mod $r$ representation associated to $X[r]$ is valued in $\GSp(2g,\F_r)$. Let $r>M$ be a prime distinct from $\ell$. Let $w$ be an extension of $v$ to the $r$-torsion field $K = F(X[r])$. By \cref{inertcrit2}, we know $X$ attains semistable reduction at $w$ over $K$. \cref{inertcrit1} now implies that the absolute inertia group at $w$ acts unipotently on $T_p(X)$ and hence also on $X[p]$. So $\rho(I_w)$ is a $p$-group. Thus the prime-to-$p$ part of $\# \rho(I_v)$ divides $\# I_v(K|F)$ which in turn divides $\# \GSp(2g,\F_r)$. Since this is true for all primes $r>M$, $r \neq \ell$, we get by \cref{Kglemma2} that the prime-to-$p$ part of $\# \rho(I_v)$ divides $K_g$.
	\end{proof}
\end{prop}

\section{Non-split Cartan subgroups of symplectic groups}
\label{sec:cartan}

Let $d \geq 1$ be an integer, and let $k$ denote the finite field of order $p^d$. Consider the symplectic bilinear form $\wedge_k$ on $k^2$ defined by
\begin{align*}
    \wedge_k(\vv_1, \vv_2) = ad-bc, \quad \text{if } \vv_1 = [a~b]^t, ~ \vv_2 = [c~d]^t.
\end{align*}
Then $\wedge = \tr_{k|\F_p} \circ \wedge_k$ is a symplectic pairing on $k^2$ valued in $\F_p$. Since the action of $\SL(2,k)$ preserves $\wedge_k$, we get an inclusion
\begin{align*}
    \SL(2,k) \rightarrow \Sp(2d,\F_p).
\end{align*}
Let $G_d \subset \GL(2,k)$ denote the subgroup of matrices with determinant in $\F_p^{\times} \subset k^{\times}$. Then $G_d$ preserves $\wedge_k$ and $\wedge$ up to scalars, and hence we get an inclusion $G_d \rightarrow \GSp(2d,\F_p)$ such that its composition with the similitude character is surjective on to $\F_p^{\times}$.

Let $l$ be the finite field of order $p^{2d}$. Then $[l:k] = 2$, and we consider an identification of $l$ with $k^2$ as vector spaces over $k$. This induces an inclusion $l^{\times} \subset \GL(2,k)$ whose image is called the non-split Cartan subgroup of $\GL(2,k)$. We consider the following subgroups
\begin{align*}
	& C = \{~ x \in l^{\times} ~|~ \Norm_{l|k}~ x \in \F_p^{\times} ~\} \subset G_d \\
	& C_1 = \{~ x \in l^{\times} ~|~ \Norm_{l|k}~ x = 1 ~\} \subset \SL(2,k).
\end{align*}
Then $C_1$ is the non-split Cartan subgroup of $\SL(2,k)$. Identifying $C$ and $C_1$ with their images under the inclusion $G_d \rightarrow \GSp(2d,\F_p)$, we see that $C \subset \GSp(2d,\F_p)$ is a cyclic subgroup of order $(p^d+1)(p-1)$ and $C_1 = C \cap \Sp(2d,\F_p)$ is the subgroup of order $p^d+1$.

The Galois group of $l$ over $\F_p$ acts naturally on $C$, with Frobenius acting by raising to $p^{th}$ power. \cref{lem:Frob2preservespairing} and \cref{lem:Frobppreservespairing} allow us to describe a subgroup $N$ inside the normalizer of $C$ in $\GSp(2d,\F_p)$, such that the action of the quotient $N/C$ on $C$ matches this Galois action. We will see that this extension is split when $p = 2$, and non-split otherwise.

\begin{lem}
\label{lem:Frob2preservespairing}
Let $p = 2$. Let $\eta \in l^{\times}$ be such that $\tr_{l|k}(\eta) = -1$, i.e., the minimal polynomial over $k$ of $\eta$ is of the form $x^2 + x + u$ for some $u \in k$. Let us identify $l$ with $k^2$ using the basis $1, \eta$.
Let $\sigma = \frob_p \in \Gal(l|\F_p)$. Then $\sigma$ acts $\F_p$-linearly on $l$, and preserves the pairing $\wedge$.
\begin{proof}
    By definition $\sigma = \frob_p$ is a $\F_p$-linear automorphism of $l$. Let $a+b\eta$ and $c+d\eta$ be two elements in $l$. With the given identification $l = k^2$, we have $\wedge(a+b\eta,c+d\eta) = \tr_{k|\F_p}(ad-bc)$. Then
    \begin{align*}
        \wedge(\sigma(a+b\eta), \sigma(c+d\eta)) & = \wedge(a^2 + b^2 (\eta^2), c^2 + d^2(\eta^2))\\
        & = \wedge((a^2-u b^2)-b^2 \eta, (c^2-u d^2)-d^2 \eta)\\
        & = \tr_{k|\F_p}(-a^2d^2+b^2c^2)\\
        & = \tr_{k|\F_p}(\frob_p(ad-bc))\\
        & = \tr_{k|\F_p}(ad-bc) = \wedge(a+b\eta,c+d\eta)
    \end{align*}
    showing that the action of $\sigma$ preserves the pairing $\wedge$.
\end{proof}
\end{lem}

It is clear that $\sigma$ has order $2d$, and the conjugation action of $\sigma$ sends $x \in C$ to $\sigma x \sigma^{-1} = \sigma(x) \cdot \sigma \circ \sigma^{-1} = \sigma(x) = x^2$. Let $N$ denote the subgroup of $\GSp(2d,\F_2)$ generated by $C$ and $\sigma$. Then $N$ is contained in the normalizer of $C$, and is isomorphic to the semidirect product $C \rtimes \langle \sigma \rangle$. In other words, we have a split short exact sequence
\begin{equation}
\label{eq:ses_split}
\begin{tikzcd}
    0 \arrow{r} & \left[ N,N \right] = C_1 = C \arrow{r} & N \arrow{r} & N^{ab} \simeq \Z/2d \arrow[bend right=30,swap]{l} \arrow{r} & 0.
\end{tikzcd}
\end{equation}
Let $x$ denote a generator of $C$, and $y = \sigma$ so that $\langle x, y | x^{2^d+1} = y^{2d} = 1, yxy^{-1} = x^2 \rangle$ is a presentation of $N$. Then the abelianization map above sends $x^ay^b \in N$ to $b \in \Z/2d$, with the obvious splitting $\Z/2d \rightarrow N$ sending $b \mapsto y^b$.

\begin{lem}
\label{lem:Frobppreservespairing}
Let $p$ be odd. Let $\eta \in l^{\times}$ such that $\eta^2 \in k^{\times}$ is a primitive root, and let us identify $l$ with $k^2$ using the basis $1, \eta$.
Let $\alpha \in l^{\times}$ and let $\sigma = \frob_p \in \Gal(l|\F_p)$. Then $\tilde{\sigma} := \alpha \sigma$ acts $\F_p$-linearly on $l$, and it preserves the pairing $\wedge$ if and only if
\[ \Norm_{l|k}(\alpha) = \eta^{1-p}. \]
Note this means that $\alpha$ can be taken to be in $k^{\times}$ if and only if $p \equiv 1 \pmod 4$.
\begin{proof}
    It is clear that $\tilde{\sigma}$ acts $\F_p$-linearly, since both $\sigma$ and multiplication by $\alpha \in l^{\times}$ are $\F_p$-linear operations.
    With the given identification $l = k^2$, we have $\wedge(a+b\eta,c+d\eta) = \tr_{k|\F_p}(ad-bc)$. If $\alpha = \alpha_1 + \alpha_2 \eta$, then we have
    \begin{align*}
        \wedge (\alpha \sigma (a+b\eta), \alpha \sigma (c+d\eta)) & = \wedge ((\alpha_1 + \alpha_2 \eta) (a^p + b^p \eta^p), (\alpha_1 + \alpha_2 \eta) (c^p + d^p \eta^p) ) \\
        &= \tr_{k|\F_p}(\eta^{p-1}(\alpha_1^2 - \alpha_2^2 \eta^2)(a^pd^p-b^pc^p))\\
        &= \tr_{k|\F_p}(\eta^{p-1}\Norm_{l|k}(\alpha)\frob_p(ad-bc))
    \end{align*}
    This is equal to $\tr_{k|\F_p}(ad-bc)$ for all $a,b,c,d \in k$ if and only if $\eta^{p-1}\Norm_{l|k}(\alpha) = 1$, which proves the lemma.
\end{proof}
\end{lem}

Suppose we choose $\eta, \alpha$ as in \cref{lem:Frobppreservespairing}, so that $\tilde{\sigma} = \alpha \sigma$ preserves $\wedge$. Then the conjugation action of $\tilde{\sigma}$ sends $x \in C$ to $\alpha \sigma x \sigma^{-1} \alpha^{-1} = \alpha x^p \alpha^{-1} = x^p$. Unlike the case of $p=2$, the order of the element $\tilde{\sigma}$ is not obvious. We have
\[ \tilde{\sigma}^n = (\alpha \sigma)^n = \alpha \alpha^p \dots \alpha^{p^{n-1}} \sigma^n = \alpha^{\frac{p^n-1}{p-1}} \sigma^n. \]
In particular, since the order of $\sigma$ is $2d$ we have
\[ \tilde{\sigma}^{2d} = \alpha^{\frac{p^{2d}-1}{p-1}} \sigma^{2d} = (\alpha^{1+p^d})^{\frac{p^d-1}{p-1}} = \Norm_{l|k}(\alpha)^{\frac{p^d-1}{p-1}} = \eta^{-(p^d-1)} = -1 \in C. \]
Hence, the element $\tilde{\sigma} \in \Sp(2d,\F_p)$ is of order $4d$, and normalizes $C$. Let $N$ denote the subgroup of $\GSp(2d,\F_p)$ generated by $C$ and $\tilde{\sigma}$. Then $N$ is contained in the normalizer of $C$, and admits a short exact sequence
\begin{align}
\label{eq:ses}
    0 \longrightarrow [N,N] = C_1 \longrightarrow N \longrightarrow N^{ab} \simeq \Z/(p-1) \times \Z/2d \longrightarrow 0.
\end{align}
Unlike the case $p=2$, this sequence does not split.
Let $x$ denote a generator of $C$, and $y = \tilde{\sigma}$ so that $N$ has the presentation $\langle x, y | x^e = 1, y^{2d} = x^{e/2}, yxy^{-1} = x^p \rangle$ where $e = (p^d+1)(p-1)$. Then, $C_1$ is generated by $x^{p-1}$ and the abelianization map above sends $x^ay^b \in N$ to $(a,b) \in \Z/(p-1) \times \Z/2d$. The similitude character $N \rightarrow \F_p^{\times}$ corresponds to the projection on to the first factor in $N^{ab}$, followed by the isomorphism $\Z/(p-1) \simeq \F_p^{\times}$ sending $1 \mapsto \Norm_{l|k}(x)$.

\section{Embedding problem}
\label{sec:embed}

Recall the subgroup $N$ contained in the normalizer of a non-split Cartan subgroup of $\GSp(2d,\F_p)$, defined in \cref{sec:cartan}. In this section, we show the existence of a number field $K$ with $\Gal(K|\Q) \simeq N$, such that the similitude character of $N$ cuts out the subfield $\Q(\zeta_p) \subset K$.

When $p = 2$, $N$ is a semi-direct product of abelian groups as shown in \cref{sec:cartan}, and furthermore the similitude condition is trivial. Hence the existence of $K$ in this case is immediate from known results on Inverse Galois problem. For example, Shafarevich's theorem \cite{Safarevic} says that every solvable group is a Galois group over $\Q$, though it is too strong for our purpose.

For the rest of this section let $p$ be an odd prime. Shafarevich's theorem again says that $N$ is a Galois group over $\Q$ since it is solvable. But this is not enough since we need additionally that our number field have $\Q(\zeta_p)$ as the appropriate subfield. So we need to study the embedding problem
\begin{equation}
\label{embedprob}
\begin{tikzcd}
    & & & G_{\Q} \arrow{d}{\phi} \arrow[dashed]{dl}[swap]{?}& \\
    0 \arrow{r} & \left[ N,N \right] \simeq \Z/(p^d+1) \arrow{r} & N \arrow{r} & \Z/(p-1) \times \Z/2d \arrow{r} & 0
\end{tikzcd}
\end{equation}
where $\phi$ is a homomorphism such that the fixed field of the kernel of $\proj_1 \circ ~\phi$ is the cyclotomic field $\Q(\zeta_p)$.
\iffalse
Let $p$ be an odd prime, and $g \geq 2$. Let $e = (p^g+1)(p-1)$. Consider the group $N = \langle x,y | x^e = 1, y^{2g} = x^{e/2}, y x y^{-1} = x^p \rangle$ and the short exact sequence
\begin{align}
	0 \longrightarrow [N,N] = \Z/(p^g+1) \longrightarrow N \longrightarrow \Z/(p-1) \times \Z/2g \longrightarrow 0.
\end{align}
The first map sends $1 \in \Z/(p^g+1)$ to $x^{p-1}$, and the second map sends $x \in N$ to $(1,0)$, and $y \in N$ to $(0,1)$.
\fi
Suppose $F | \Q$ is a number field such that $F \cap \Q(\zeta_p) = \Q$ and $\Gal(F|\Q) \simeq \Z/2d$. Then, $F(\zeta_p) | \Q$ is Galois over $\Q$ with Galois group isomorphic to $\Z/(p-1) \times \Z/2d$. Let $\phi$ be the homomorphism cutting out $F(\zeta_p)$ i.e., $\Qbar^{\ker \phi} = F(\zeta_p)$. The embedding problem given by (\ref{embedprob}) asks whether $\phi$ can be lifted to a map $\tilde{\phi} : G_{\Q} \rightarrow N$ such that the diagram commutes. Such a lift $\tilde{\phi}$ describes an embedding of $F(\zeta_p)$ into a number field $L = \Qbar^{\ker \tilde{\phi}}$ with $\Gal(L|\Q) \subseteq N$. A lift $\tilde{\phi}$ is called a proper solution to the embedding problem if it is surjective, i.e., if $\Gal(L|\Q) \simeq N$. We refer to \cite[\S 3.5]{Neukirch2008} for a detailed discussion of embedding problems.

The rest of this section is devoted to proving the existence of a proper solution to the embedding problem given by (\ref{embedprob}) for a suitably chosen initial field $F$. We follow the general strategy to study these types of problems. Let $\epsilon$ denote the cohomology class in $H^2(\Z/(p-1) \times \Z/2d, \Z/(p^d+1))$ corresponding to the group extension $N$ in \cref{eq:ses}. Then there exists a lift $\tilde{\phi}$ if and only if $\phi^* \epsilon = 0 \in H^2(\Q,\Z/(p^d+1))$ \cite[Prop 3.5.9]{Neukirch2008}. We show $\phi^* \epsilon = 0$ in two steps. First, we show that the restriction $\res_{\ell} (\phi^* \epsilon) = 0 \in H^2(\Q_{\ell}, \Z/(p^d+1))$ for all primes $\ell$ including the infinite prime. Second, we show that Hasse principle holds in our case. That is, if all the local restrictions of a global cohomology class are trivial, then the global cohomology class itself is trivial. Finally, we exploit the fact \cite[Prop 3.5.11]{Neukirch2008} that the space of solutions to the embedding problem given by (\ref{embedprob}) is a principal homogenous space over $H^1(\Q,\Z/(p^d+1))$, and twist using a suitable class to obtain properness.

We will choose $F$ so that all ramification in $F$ is tame and all the local embedding problems are solvable. Let $2d = 2^n d_1$ where $d_1$ is odd. Then $\Z/2d \simeq \Z/2^n \times \Z/d_1$. We will choose Galois extensions $F_1$ and $F_2$ of $\Q$ with Galois groups $\Z/2^n$ and $\Z/d_1$ respectively, and define $F$ to be their compositum. For $i = 1,2$, we take $F_i$ to be the unique subfield of the above mentioned degree inside the cyclotomic field $\Q(\zeta_{N_i})$, for certain primes $N_i$ described below.

Let $N_2$ be a prime such that $N_2 \equiv 1 \pmod {d_1}$, say $N_2 = 2 \alpha d_1 + 1$ for some $\alpha \in \N$. Let $N_1$ be a prime satisfying
\begin{enumerate}[(a)]
    \item \label{conda} $N_1 \equiv 2^n+1 \pmod {2^{n+1}}$.
    \item \label{condb} $N_1 \equiv 1 \pmod {N_2}$.
    \item \label{condc} $p \not\equiv \square \pmod {N_1}$.
\end{enumerate}
The third condition can be rewritten as a congruence condition on $N_1$ modulo $p$ using quadratic reciprocity. Dirichlet's theorem on primes in arithmetic progression guarantees the existence of such primes $N_1,N_2$.

Let $S$ denote the finite set containing the infinite prime and all primes ramified in $F(\zeta_p) | \Q$, i.e., $S = \Ram(F|\Q) \cup \{ \infty, p \}$. With the choices made above, we have $\Ram(F|\Q) = \{N_1, N_2\}$, and $S = \{ \infty, p, N_1, N_2 \}$ and $F(\zeta_p)$ is tamely ramified at each finite prime in $S$. For primes not in $S$, the extension $F(\zeta_p) | \Q$ is unramified. So the local embedding problems at these primes are trivially solvable, meaning that there is no local obstruction to the embedding problem given by (\ref{embedprob}) at these primes.
We now study the local embedding problem for each prime in $S$.

\subsection{Local obstruction at $\infty$}
\label{subsec:infty}

If $F_1$ is chosen as above, condition \ref{conda} on $N_1$ implies that $F_1$ is not a totally real extension of $\Q$. That is, complex conjugation is given by the non-trivial order $2$ element in $\Gal(F_1|\Q)$. Complex conjugation acts trivially on $F_2$ since the order of $\Gal(F_2|\Q) = \deg(F_2) = d_1$ is odd. Thus, complex conjugation in $\Gal(F(\zeta_p)|\Q) = \Z/(p-1) \times \Z/2d$ is given by the element $(\frac{p-1}{2}, d)$.

The element $x^{(p-1)/2}y^d$ is clearly a lift of complex conjugation to $N$. Recalling that $e = (p^d+1)(p-1)$, we further have
\begin{align*}
	\left(x^{\frac{p-1}{2}}y^d\right)^2 = x^{\frac{p-1}{2}} \left(y^d x^{\frac{p-1}{2}} y^{-d}\right) y^{2d} = x^{\frac{p-1}{2}} x^{\frac{(p-1)p^d}{2}} y^{2d} = x^{\frac{e}{2}} y^{2d} = 1,
\end{align*}
so the lift has order $2$. This shows that there is no local obstruction at the infinite place to the embedding problem given by (\ref{embedprob}).

\subsection{Local obstruction at $p$}
\label{subsec:p}

The local obstruction at $p$ is measured by whether or not the restriction of $\phi$ to the decomposition group $G_{\Q_p}$, can be lifted to a homomorphism $G_{\Q_p} \rightarrow N$. The map $\phi|_{G_{\Q_p}}$ factors through the tame Galois group $G_{\Q_p}^{\text{tame}}$ which is a profinite group with presentation $\langle \sigma,\tau | \sigma \tau \sigma^{-1} = \tau^p \rangle$, where $\tau$ is a generator of tame inertia, and $\sigma$ is a lift of the Frobenius of the maximal unramified extension. Without loss of generality, suppose that $\phi$ sends $\sigma$ to $(0,a)$ and $\tau$ to $(1,0)$ in $\Gal(F(\zeta_p)|\Q) \simeq \Z/(p-1) \times \Z/2d$.

\begin{prop}
\label{prop:localp}
There exist $\tilde{\sigma}, \tilde{\tau} \in N$ lifting $(0,a)$ and $(1,0)$ and satisfying $\tilde{\sigma}\tilde{\tau}\tilde{\sigma}^{-1} = \tilde{\tau}^p$ if and only if $a \equiv 1 \pmod 2$.
\end{prop}
\begin{proof}
    Let $\tilde{\sigma} = x^{l(p-1)}y^a$ and $\tilde{\tau} = x^{1+k(p-1)}$ be any lifts. We have
    \begin{align*}
        \tilde{\sigma}\tilde{\tau}\tilde{\sigma}^{-1}\tilde{\tau}^{-p} = y^a x^{1+k(p-1)} y^{-a} x^{-(1+k(p-1))p} = x^{(1+k(p-1))(p^a-p)}
    \end{align*}
    If $a = 0$, then the desired condition $\tilde{\sigma}\tilde{\tau}\tilde{\sigma}^{-1} = \tilde{\tau}^p$ cannot be met since the equation
    \[ 1 + k(p-1) \equiv 0 \pmod {p^d+1} \]
    has no solution.
    If $a = 1$, any choice of $k$ and $l$ gives desired lifts.
    
    Assume $a \geq 2$. We get a lift satisfying $\tilde{\sigma}\tilde{\tau}\tilde{\sigma}^{-1} = \tilde{\tau}^p$ if and only if there exists $k \in \Z/(p^d+1)$ satisfying the equation
    \begin{align*}
        (1+k(p-1))(p^a-p) & \equiv 0 \pmod e \\
        \text{i.e., } k(p-1) & \equiv -1 \pmod {e'}
    \end{align*}
    where
    \[ e' = \frac{e}{\gcd(e,p^a-p)} = \frac{p^d+1}{\gcd(p^d+1,1+p+p^2+\dots+p^{a-2})}.\]
    This equation has a solution if and only if $p-1$ is invertible modulo $e'$. Since $\gcd(p-1,e')$ divides $\gcd(p-1,p^d+1) = 2$, 
    this happens if and only if $e'$ is odd.

    If $d$ is even, then $p^d+1 \equiv 2 \pmod{4}$. Hence $e'$ is odd if and only if $2$ divides $1 + p + p^2 + \dots + p^{a-2}$, which happens if and only if $a \equiv 1 \pmod 2$.
    
    Suppose $d$ is odd. Let $m \geq 1$ be such that $p \equiv 2^m-1 \pmod {2^{m+1}}$. Then
    $p^d + 1 \equiv p + 1 \equiv 2^m \pmod {2^{m+1}}$.
    Hence, $e'$ is odd if and only if $2^m$ divides $1 + p + p^2 + \dots + p^{a-2}$. Since $p \equiv -1 \pmod{2^m}$, this happens if and only if $a \equiv 1 \pmod 2$.
    This completes the proof of the proposition.
\end{proof}

The proposition says that the local obstruction at $p$ to the embedding problem vanishes if and only if $\frob_p$ is not a square in $\Gal(F|\Q)$. Equivalently $\frob_p$ is not a square in $\Gal(F_1|\Q)$. This holds as a result of condition \ref{condc}.

\subsection{Local obstruction at $N_1$}
\label{subsec:N1}

The prime $N_1$ is unramified in $\Q(\zeta_p)$, totally tamely ramified in $F_1 \subset \Q(\zeta_{N_1})$, and split in $F_2$ thanks to condition \ref{condb}. Hence the restriction of $\phi$ to the decomposition group ${G_{\Q_{N_1}}}$ factors through the profinite tame quotient $G_{\Q_{N_1}}^{\text{tame}} = \langle \sigma,\tau | \sigma \tau \sigma^{-1} = \tau^{N_1} \rangle$ as before, and without loss of generality, we may suppose that $\phi$ sends $\sigma$ to $\frob_{N_1} = (a,0)$ and $\tau$ to $(0,d_1)$ in $\Gal(F(\zeta_p)|\Q) \simeq \Z/(p-1) \times \Z/2d$.

Note that conditions \ref{conda} and \ref{condc} already determine the parity of $a$ by quadratic reciprocity. To be precise, if $d$ is even making $n \geq 2$ and hence $N_1 \equiv 1 \pmod 4$, or if $p \equiv 1 \pmod 4$, then by quadratic reciprocity we have that $N_1 \not \equiv \square \pmod p$ meaning that $a$ is odd. If $d$ is odd and $p \equiv 3 \pmod 4$ then $a$ is even. This will be used in the following proposition which proves the existence of elements $\tilde{\sigma}, \tilde{\tau} \in N$ satisfying certain properties. These elements determine a homomorphism $G_{\Q_{N_1}} \rightarrow N$ factoring through the tame Galois group that lifts $\phi|_{G_{\Q_{N_1}}}$. This shows that there is no local obstruction to the embedding problem at the prime $N_1$.

\begin{prop}
\label{prop:localN1}
There exist $\tilde{\sigma}, \tilde{\tau} \in N$ satisfying $\tilde{\sigma}\tilde{\tau}\tilde{\sigma}^{-1} = \tilde{\tau}^{N_1}$ and lifting the elements $\phi(\sigma) = (a,0)$ and $\phi(\tau) = (0,d_1)$ in $\Z/(p-1) \times \Z/2d$.
\end{prop}
\begin{proof}
Consider the elements $\tilde{\sigma} = x^{a+k(p-1)}$ and $\tilde{\tau} = y^{d_1}$ in the group $N$ lifting the elements $\phi(\sigma)$ and $\phi(\tau)$. We will show that there is a choice of $k$ so that
$\tilde{\sigma} \tilde{\tau} \tilde{\sigma}^{-1} = \tilde{\tau}^{N_1}$. 
We first simplify both sides of the expression.
\begin{align*}
    \tilde{\sigma} \tilde{\tau} \tilde{\sigma}^{-1} = x^{a+k(p-1)} \left( y^{d_1} x^{-(a+k(p-1))} y^{-d_1} \right) y^{d_1} = x^{(a+k(p-1))(1-p^{d_1})} y^{d_1}.
\end{align*}
Since the order of $y \in N$ is $4d = 2^{n+1}d_1$, and condition \ref{conda} says that $N_1 \equiv 2^n+1 \pmod{2^{n+1}}$,
\begin{align*}
    \tilde{\tau}^{N_1} = y^{N_1 d_1} = y^{(2^n+1)d_1} = y^{2d}y^{d_1} = x^{e/2}y^{d_1}.
\end{align*}
Thus, we need to show that there is a solution $k$ to the equation
\begin{align*}
    (a+k(p-1))(1-p^{d_1}) & \equiv e/2 \pmod{e}\\
    \text{i.e., } k(p-1)(1-p^{d_1}) & \equiv e/2 - a(1-p^{d_1}) \pmod{e}\\
    \text{i.e., } k(1-p^{d_1}) & \equiv \frac{p^d+1}{2} + a(1+p+p^2+\dots+p^{d_1-1}) \pmod{p^d+1}.
\end{align*}
Since $d_1$ divides $d$, it is clear that $\gcd(1-p^{d_1}, p^d+1) = 2$. Hence, the above equation has a solution if and only if
\begin{align*}
	\frac{p^d+1}{2} + a(1+p+p^2+\dots+p^{d_1-1}) & \equiv 0 \pmod{2}\\
	i.e., \quad \frac{p^d+1}{2} + a & \equiv 0 \pmod{2} \quad (\text{since } d_1 \text{ is odd})
\end{align*}
The parity condition on $a$ we described earlier ensures that this holds. If $d$ is even or $p \equiv 1 \pmod 4$, then both $\frac{p^d+1}{2}$ and $a$ are odd. Otherwise, both are even.
\end{proof}

\subsection{Local obstruction at $N_2$}
\label{subsec:N2}

The prime $N_2$ is unramified in $\Q(\zeta_p)$ and $F_1$, and totally tamely ramified in $F_2 \subset \Q(\zeta_{N_2})$. Hence the restriction of $\phi$ to the decomposition group $G_{\Q_{N_2}}$ factors through the profinite tame quotient $G_{\Q_{N_2}}^{\text{tame}} = \langle \sigma,\tau | \sigma \tau \sigma^{-1} = \tau^{N_2} \rangle$, and without loss of generality, we may suppose that $\phi$ sends $\sigma$ to $\frob_{N_2} = (a,bd_1)$ and $\tau$ to $(0,2^n)$ in $\Gal(F(\zeta_p)|\Q) \simeq \Z/(p-1) \times \Z/2d$.

\begin{prop}
\label{prop:localN2}
There exist $\tilde{\sigma}, \tilde{\tau} \in N$ satisfying $\tilde{\sigma}\tilde{\tau}\tilde{\sigma}^{-1} = \tilde{\tau}^{N_2}$ and lifting the elements $\phi(\sigma) = (a,bd_1)$ and $\phi(\tau) = (0,2^n)$ in $\Z/(p-1) \times \Z/2d$.
\end{prop}
\begin{proof}
Consider the elements $\tilde{\sigma} = x^{a+k(p-1)}y^{bd_1}$ and $\tilde{\tau} = y^{2^n}$ in the group $N$ lifting the elements $\phi(\sigma)$ and $\phi(\tau)$. We will show that there is a choice of $k$ so that $\tilde{\sigma} \tilde{\tau} \tilde{\sigma}^{-1} = \tilde{\tau}^{N_2}$.
We first simplify both sides of the expression, recalling that $N_2 = 2 \alpha d_1 + 1$.
\begin{align*}
    \tilde{\sigma} \tilde{\tau} \tilde{\sigma}^{-1} &= x^{a+k(p-1)} \left( y^{2^n} x^{-(a+k(p-1))} y^{-2^n} \right) y^{2^n} = x^{(a+k(p-1))(1-p^{2^n})} y^{2^n}.\\
    \tilde{\tau}^{N_2} &= y^{2^n N_2} = y^{2^{n+1}d_1 \alpha} y^{2^n} = y^{4d\alpha}y^{2^n} = y^{2^n}.
\end{align*}
Thus we need to show that there is a solution $k$ to the equation
\begin{align*}
    (a+k(p-1))(1-p^{2^n}) & \equiv 0 \pmod{e}.\\
    \text{i.e., } k(1-p^{2^n}) & \equiv a(1+p+p^2+\dots+p^{2^n-1}) \pmod{p^d+1}.\\
    \text{i.e., } k(1-p^{2^{n-1}})(1+p^{2^{n-1}}) & \equiv a(1+p^{2^{n-1}})(1+p+\dots+p^{2^{n-1}-1}) \pmod{p^d + 1}.
\end{align*}
Since $d = 2^{n-1} d_1$ and $d_1$ is odd, this equation is the same as
\begin{align*}
    k(1-p^{2^{n-1}}) & \equiv a(1+p+\dots+p^{2^{n-1}-1}) \pmod{M},
\end{align*}
where
\begin{align*}
    M = \frac{p^d+1}{p^{2^{n-1}}+1} = \sum_{i = 0}^{d_1-1} \left(- p^{2^{n-1}}\right)^i.
\end{align*}
Since $M \equiv \sum_{i = 0}^{d_1-1} (-1)^i \equiv 1 \pmod{1-p^{2^{n-1}}}$, we know that $1-p^{2^{n-1}}$ is coprime to $M$. Hence the above equation does indeed have a solution.
\end{proof}

The elements $\tilde{\sigma}, \tilde{\tau}$ from \cref{prop:localN2} determine a homomorphism $G_{\Q_{N_2}} \rightarrow N$ factoring through the tame Galois group, that lifts $\phi|_{G_{\Q_{N_2}}}$.
Hence there is no local obstruction to the embedding problem at the prime $N_2$ either.

\subsection{Global obstruction}

Let $A$ denote the $G_{\Q}$-module $[N,N] = \Z/(p^d+1)$. Recall that the Galois action factors through the map $\phi : G_{\Q} \rightarrow N^{ab} \simeq \Z/(p-1) \times \Z/2d$ and is given by conjugation in $N$ as in the short exact sequence \cref{eq:ses}. Note that this action further factors through $\proj_2 \circ ~ \phi : G_{\Q} \rightarrow \Gal(F|\Q) \simeq \Z/2d$ and a generator of $\Z/2d$ acts on $A$ as raising to $p^{th}$ power. Global obstruction to the embedding problem given by (\ref{embedprob}) is measured by the group $\Sh_{\Q}^2(A)$ defined as
\[ \Sh_{\Q}^{2}(A) = \ker \left( H^2(\Q,A) \longrightarrow \prod_v H^2(\Q_v,A) \right), \]
where $v$ runs over all places of $\Q$.

\begin{thm}
    \label{thm:no-global-obs}
    There is no global obstruction to this embedding problem, i.e., $\Sh_{\Q}^2(A) = 0$.
\end{thm}
\begin{proof}
    By Poitou-Tate duality, we have $\Sh_{\Q}^2(A) \simeq \Sh_{\Q}^1(A^{\vee})^{\vee}$, where $A^{\vee} = \Hom(A,\overline{\Q}^{\times})$ is the dual module. If we let $m = p^d+1$, then 
    $A^{\vee} = \Hom(A,\mu_m)$. Let $k$ be the trivializing extension of $A^{\vee}$.
    Consider the following commutative diagram
    \begin{equation}
        \label{commdiag1}
        \begin{tikzcd}
            & & H^1(k,A^{\vee}) \arrow{r} & \prod\limits_w H^1(k_w,A^{\vee})\\
            0 \arrow{r} & \Sh_{\Q}^1(A^{\vee}) \arrow{r} & H^1(\Q,A^{\vee}) \arrow[u] \arrow{r} & \prod\limits_v H^1(\Q_v,A^{\vee}) \arrow[u] \\
            0 \arrow{r} & \Sh_{k|\Q}^1(A^{\vee}) \arrow[u] \arrow{r} & H^1(k|\Q, A^{\vee}) \arrow[u] \arrow{r} & \prod\limits_v H^1(k_v|\Q_v,A^{\vee}) \arrow[u]\\
            & & 0 \arrow[u] & 0 \arrow[u]
        \end{tikzcd}
    \end{equation}
    where the vertical maps are coming from the inflation restriction sequence. Since $A^{\vee}$ is trivial as a $G_k$-module, and $w$ ranges over all places of $k$, Hasse principle holds for the $G_k$-module $A^{\vee}$ as per \cite[Thm 9.1.9.(i)]{Neukirch2008}. That is, the horizontal map at the top is injective. Thus we get the isomorphism
    \begin{align}
        \label{coh-of-fin-grp}
        \Sh_{\Q}^1(A^{\vee}) \simeq \Sh_{k|\Q}^1(A^{\vee}),
    \end{align}
    bringing us to the study of the cohomology of the module $A^{\vee}$ of the finite group $\Gal(k|\Q)$.
    
    Since the Galois action on $A$ factors through $F$, we know that $k$ is contained in $F(\zeta_m)$. In fact it is easy to see that $A^{\vee}$ as a $\Gal(k|\Q)$-module is isomorphic to $\mu_m$ as a $\Gal(\Q(\zeta_m)|\Q)$-module. That is, there is an isomorphism of pairs
    \begin{align}
        \label{transfer}
        \psi : \left( \Gal(k|\Q), ~ A^{\vee} \right) \longrightarrow \left( \Gal(\Q(\zeta_m)|\Q), \mu_m \right).
    \end{align}
    The map $\Gal(F(\zeta_m)|\Q) \simeq \Z/2d \times \left( \Z/m \right)^{\times} \rightarrow \left( \Z/m \right)^{\times} \simeq \Gal(\Q(\zeta_m)|\Q)$ sending $(a,b) \mapsto p^{-a}b$ induces the isomorphism $\psi$ on the groups. We note here that since inertia subgroup behaves well with respect to quotients, for any prime unramified in $F|\Q$, the isomorphism $\psi$ identifies the inertia subgroups of $k$ and $\Q(\zeta_m)$ at that prime. In particular, the inertia subgroups at $2$ get identified.
    
    In order to study $\Sh_{k|\Q}^1(A^{\vee})$, we will use the isomorphism $\psi$ in \cref{transfer} and instead study the familiar module $\mu_m$.
    Before doing this, we relax local conditions slightly.
    Let $T$ denote the set of all odd primes that are unramified in $k|\Q$. Let $\mathcal{L}$ denote the Selmer conditions given by
    \begin{align*}
        L_v = \begin{cases} H^1_{ur}(k_v|\Q_v, A^{\vee}), \quad & \text{if } v = 2\\
        0, \quad & \text{if } v \in T \\
        H^1(k_v|\Q_v,A^{\vee}), \quad & \text{otherwise} \end{cases}
    \end{align*}
    In words, the local condition at $2$ is relaxed from split to unramified, and the local conditions at ramified primes are fully relaxed. The corresponding Selmer group $H_{\mathcal{L}}^1(k|\Q,A^{\vee})$ is defined as
    \begin{align*}
        H_{\mathcal{L}}^1(k|\Q,A^{\vee}) = \ker \left( H^1(k|\Q,A^{\vee}) \longrightarrow \prod_v H^1(k_v|\Q_v,A^{\vee})/L_v \right)
    \end{align*}
    and it clearly contains $\Sh_{k|\Q}^1(A^{\vee})$. So, it is enough to show that $H_{\mathcal{L}}^1(k|\Q,A^{\vee}) = 0$.
    
    The Selmer condition $\mathcal{L}$ amounts exactly to requiring that restriction to inertia subgroup at $2$ of $k|\Q$, and to any cyclic subgroup of $\Gal(k|\Q)$ is zero.
    As mentioned earlier, the isomorphism $\psi$ in \cref{transfer} identifies the inertia group at $2$ of $k|\Q$ with that of $\Q(\zeta_m)|\Q$. Thus, the induced isomorphism in group cohomology $\psi^* : H^1(k|\Q,A^{\vee}) \simeq H^1(\Q(\zeta_m)|\Q, \mu_m)$ gives an isomorphism of Selmer groups
    \begin{align}
        \label{transferselmer}
        H_{\mathcal{L}}^1(k|\Q,A^{\vee}) \simeq H_{\mathcal{L}'}^1(\Q(\zeta_m)|\Q, \mu_m),
    \end{align}
    where $\mathcal{L}'$ is a similar set of Selmer conditions. To be precise, let $T'$ denote the set of all odd primes that are unramified in $\Q(\zeta_m)|\Q$. Then, $\mathcal{L}'$ imposes the unramified condition at the prime $2$, and the split condition at every prime in $T'$.
    
    We temporarily forget the condition at $2$, and consider a commutative diagram similar to (\ref{commdiag1}) for the Galois module $\mu_m$ and the set $T'$
    \begin{equation}
        \label{commdiag2}
        \begin{tikzcd}
            & & H^1(\Q(\zeta_m),\mu_m) \arrow{r} & \prod\limits_{w \in T' } H^1(\Q(\zeta_m)_w,\mu_m)\\
            0 \arrow{r} & \Sh_{\Q}^1(T',\mu_m) \arrow{r} & H^1(\Q,\mu_m) \arrow[u] \arrow{r} & \prod\limits_{v \in T'} H^1(\Q_v,\mu_m) \arrow[u] \\
            0 \arrow{r} & \Sh_{\Q(\zeta_m)|\Q}^1(T',\mu_m) \arrow[u] \arrow{r} & H^1(\Q(\zeta_m)|\Q, \mu_m) \arrow[u] \arrow{r} & \prod\limits_{v \in T'} H^1(\Q(\zeta_m)_v|\Q_v,\mu_m) \arrow[u]\\
            & & 0 \arrow[u] & 0 \arrow[u]
        \end{tikzcd}
    \end{equation}
    Then \cite[Thm 9.1.9]{Neukirch2008} again says that the horizontal map at the top is injective, and hence
    \begin{align}
        \label{selmerinf}
        \Sh_{\Q(\zeta_m)|\Q}^1(T',\mu_m) \simeq \Sh_{\Q}^1(T',\mu_m).
    \end{align}
    Furthermore, the same theorem says that Hasse principle for $\mu_m$ holds over $\Q$ as long as we are not in a special case. In fact, the obstruction to Hasse principle is described precisely.
    \begin{align*}
        \Sh_{\Q}^1(T',\mu_m) = \begin{cases}
            0, \quad & \text{if } (\Q,m,T') \text{ is not a special case}\\
            \Z/2, \quad & \text{if } (\Q,m,T') \text{ is a special case}
        \end{cases}
    \end{align*}
    As per the remarks following \cite[Lemma 9.1.8]{Neukirch2008}, the special case is equivalent to the statement that $8$ divides $m$, since $T'$ only consists of odd primes and has Dirichlet density $1$.
    
    If $(\Q,m,T')$ is not a special case, then we are done by \cref{transferselmer}, \cref{selmerinf} and the inclusion
    \begin{align*}
        H_{\mathcal{L'}}^1(\Q(\zeta_m)|\Q, \mu_m) \subseteq \Sh_{\Q(\zeta_m)|\Q}^1(T',\mu_m).
    \end{align*}
    Suppose $(\Q,m,T')$ is a special case. Then $8$ divides $m$. Let $m = 2^r m_1$ with $m_1$ odd and $r \geq 3$. Then the non-trivial element in $\Sh_{\Q}^1(T',\mu_m) \simeq \Z/2$ is the inflation of the class in $H^1(\Q(\zeta_m)|\Q, \mu_m)$ represented by the cocycle
    \[ \Gal(\Q(\zeta_m)|\Q) \longrightarrow \Gal(\Q(\zeta_{2^r})|\Q) \longrightarrow \Gal(\Q(\sqrt{-2})|\Q) \xrightarrow{\simeq} \{ \pm 1 \} \subseteq \mu_m. \]
    It is non-trivial when restricted to $\Gal(\Q(\zeta_m)|\Q(\zeta_{m_1}))$, which is the inertia group at $2$ of $\Q(\zeta_m)|\Q$. Hence, this class fails the unramified condition at $2$ of $\mathcal{L}'$. So, we get that $H_{\mathcal{L}'}^1(\Q(\zeta_m)|\Q, \mu_m) = 0$, and we are done by \cref{transferselmer}.
\end{proof}

We have now shown that the homomorphism $\phi : G_{\Q} \rightarrow \Gal(F(\zeta_p)|\Q) \simeq \Z/(p-1) \times \Z/2d$ in (\ref{embedprob}) lifts to some homomorphism $\tilde{\phi} : G_{\Q} \rightarrow N$. The map $\tilde{\phi}$ is not necessarily surjective. But we can twist it using a suitable cohomology class in $H^1(\Q,A)$ to get a surjective lift, as proved below. If $c: G_{\Q} \rightarrow A \subset N$ is a representing cocycle, then the twisted solution it determines is given by $c \cdot \tilde{\phi}$.

Choose a prime $v$ that splits completely in both $k$ and $F(\zeta_p)$. Cebotarev density theorem guarantees the existence of such a prime. The fact that $v$ splits completely in $F(\zeta_p)$ implies that $\phi$ is trivial on $G_{\Q_v}$ and hence the restriction of $\tilde{\phi}$ to ${G_{\Q_v}}$ lands inside $A$. So, the map $\tilde{\phi}|_{G_{\Q_v}} \in \Hom(G_{\Q_v},A) = H^1(\Q_v,A)$. Choose another homomorphism $c_v \in H^1(\Q_v,A)$ so that $c_v \cdot \tilde{\phi} : G_{\Q_v} \rightarrow A$ is surjective. If there is a global cohomology class in $H^1(\Q,A)$ which restricts to $c_v \in H^1(\Q_v,A)$, then twisting by this class gives us a proper solution. The existence of such a class is guaranteed by the following proposition.

\begin{prop}
    \label{properness}
    The map $H^1(\Q,A) \longrightarrow H^1(\Q_v,A)$ is surjective.
\end{prop}
\begin{proof}
    Let $\coker_{\Q}^1(T,M)$ denote the cokernel of the restriction map
    \begin{align*}
        H^1(\Q,M) \rightarrow \prod\limits_{v \in T} H^1(\Q_v,M).
    \end{align*}
    for a $G_{\Q}$-module $M$ and a set $T$ of places of $\Q$. We want to show $\coker_{\Q}^1(\{v\},A) = 0$. According to \cite[Lemma 9.2.2]{Neukirch2008}, there is a canonical short exact sequence
    \begin{align*}
        0 \longrightarrow \Sh_{\Q}^1(A^{\vee}) \longrightarrow \Sh_{\Q}^1(S \setminus \{v\}, A^{\vee}) \longrightarrow \coker_{\Q}^1(\{v\},A)^{\vee} \longrightarrow 0
    \end{align*}
    where $S$ is the set of all places of $\Q$. So it is enough to show that 
    $\Sh_{\Q}^1(S \setminus \{v\}, A^{\vee}) = \Sh_{\Q}^1(A^{\vee})$.
    Following a similar argument as in the proof of \cref{thm:no-global-obs}, we get that
    \[ \Sh_{\Q}^1(S \setminus \{v\}, A^{\vee}) \simeq \Sh_{k|\Q}^1(S \setminus \{v\}, A^{\vee}) \]
    Since the prime $v$ was chosen to be split in $k|\Q$, the decomposition group of $v$ inside $\Gal(k|\Q)$ is trivial. Hence the restriction map at $v$ on the finite group cohomology $H^1(\Gal(k|\Q),A^{\vee})$ is automatically zero, meaning that the local condition at $v$ is vacuous. Thus we have the isomorphism
    \[ \Sh_{\Q}^1(S \setminus \{v\}, A^{\vee}) \simeq \Sh_{k|\Q}^1(S \setminus \{v\}, A^{\vee}) = \Sh_{k|\Q}^1(A^{\vee}) %= 0
    \simeq \Sh_{\Q}^1(A^{\vee}) \]
    from which we deduce that $\coker_{\Q}^1(\{v\},A) = 0$.
\end{proof}

\section{Proof of main theorem}
\label{sec:proof-main-theorem}

We first discuss some preliminaries about the desired condition on order of inertia $\rho(I_{\ell})$ at an auxiliary prime $\ell$. 
Recall the subgroups $C, C_1$ and $N$ of $\GSp(2d,\F_p)$ introduced in \cref{sec:cartan}. For every positive integer $d \leq g$, since $\GSp(2d,\F_p) \subset \GSp(2g,\F_p)$, the group $\GSp(2g,\F_p)$ contains cyclic subgroups $C, C_1$ of orders $(p^d+1)(p-1)$ and $p^d+1$, and a subgroup $N$ of the normalizer of $C$ as described in \cref{sec:cartan} such that $[N,N] = C_1$.
We desire $\rho(I_{\ell}) \subset [N,N]$ to have a prime power order $q$ not dividing $K_g$. So we first look for prime powers $q$ that divide $p^d+1$ for some $1 \leq d \leq g$, but do not divide $K_g$.

\begin{lem}
\label{Zsig}
    Let $g \geq 7$ and $p$ be any prime. Then, there exists a prime $q > 2g+1$ such that $q$ divides $p^d+1$ for some $1 \leq d \leq g$.
\end{lem}
\begin{proof}
    Let $\pi$ denote the prime counting function. Zsigmondy's theorem implies that for any prime $p$ and $n \geq 1$, with the exception of $p = 2, n = 3$, there is a prime divisor of $p^n+1$ which does not divide $p^m+1$ for any $m < n$ \cite{Zsigmondy}.\\
    \textbf{Case 1: $p \neq 2$}.\\
    If $g \geq 7$ then $\pi(2g+1) \leq g-1$. Zsigmondy's theorem implies that there are at least $g$ distinct prime numbers that divide some number in the set $\{p^d+1 : 1 \leq d \leq g\}$. So one of them has to be bigger than $2g+1$.\\
    \textbf{Case 2: $p = 2$}.\\
    If $g = 7, 8, 9$, we may take $d = 7$ and $q = 43$. If $g \geq 10$ then $\pi(2g+1) \leq g-2$. Zsigmondy's theorem implies that there are at least $g-1$ distinct prime numbers that divide some number in the set $\{2^d+1 : 1 \leq d \leq g\}$. So one of them has to be bigger than $2g+1$.
\end{proof}

\cref{Kglemma1} and \cref{Zsig} ensure that when $g \geq 7$ there exists a prime $q$ that divides $p^d+1$ for some $1 \leq d \leq g$, and does not divide $K_g$.
Suppose $2 \leq g \leq 6$. If $p$ is a large enough prime, for example, if $p^g+1 > K_g$, then there exists a prime power $q$ that divides $p^g+1$ and does not divide $K_g$. This leaves only finitely many cases $(g,p)$ to be dealt with. For each of them except $(g,p) = (2,2),$ $(2,3),$ $(3,2),$ $(3,3)$ we check explicitly that there exists some prime power $q$ dividing $p^d+1$ for some $1 \leq d \leq g$, such that $q$ does not divide $K_g$.

We are now ready to prove \cref{thm:main-theorem}.
\begin{proof}
    Suppose $(g,p) \neq (3,3)$. By the preceding discussion, we find a number $d$ and a prime power $q$ such that $1 \leq d \leq g$, $q$ divides $p^d+1$, and $q$ does not divide $K_g$. Let $C, C_1, N$ denote the subgroups of $\GSp(2d,\F_p)$ of orders $(p^d+1)(p-1)$, $p^d+1$ and $2d(p^d+1)(p-1)$ as defined in \cref{sec:cartan}. We will consider them as subgroups of $\GSp(2g,\F_p)$ by a fixed inclusion $\GSp(2d,\F_p) \subset \GSp(2g,\F_p)$ Choose a number field $F$ and define $k$ to be the trivializing extension of the dual module $A^{\vee}$, just as in \cref{sec:embed}. The calculations in \cref{sec:embed} say that there is no obstruction to the embedding problem given by (\ref{embedprob}). Let $\tilde{\phi} : G_{\Q} \rightarrow N$ be a solution to this embedding problem.
    
    In order to get desired inertia at an auxiliary prime, we follow the same approach that was used in \cref{sec:embed} to get properness. In addition to the prime $v$ and the cohomology class $c_v \in H^1(\Q_v,A)$ chosen in \cref{sec:embed} to get properness, choose an auxiliary prime $\ell \equiv 1 \pmod q$ that splits completely in $k$ and $F(\zeta_p)$, and a homomorphism $c_{\ell} : G_{\Q_{\ell}} \rightarrow A=[N,N]$ so that the image of 
    $I_{\ell}$ under $c_l \cdot \tilde{\phi}$ is the cyclic subgroup of $[N,N]$ of order $q$. The proof of \cref{properness} goes through to show that the restriction map
    \[ H^1(\Q,A) \rightarrow H^1(\Q_v,A) \times H^1(\Q_{\ell},A) \]
    is surjective. Thus, there is a global cohomology class $c \in H^1(\Q,A)$ which restricts to $c_v$ and $c_l$. Twisting $\tilde{\phi}$ by this class produces a representation $\rho: G_{\Q} \twoheadrightarrow N \subset \GSp(2d,\F_p) \subset \GSp(2g,\F_p)$ with $\#\rho(I_{\ell}) = q \nmid K_g$. \cref{inertialorderdividesKg} now implies that $\rho$ does not arise as the $p$-torsion representation of any abelian variety over $\Q$.
    
    Suppose $(g,p) = (3,3)$. We find that there is a subgroup $N \subset \GSp(6,\F_3)$ of order 78, with surjective similitude character. It is a semi-direct product of $\Z/13$ and $\Z/6$ with presentation $\langle x,y | x^{13} = y^6 = 1, yxy^{-1} = x^4 \rangle$. We take $\phi : G_{\Q} \rightarrow \Gal(\Q(\zeta_9)|\Q) \simeq \Z/6$. Since $N$ is a semi-direct product, the resulting embedding problem is trivially solvable. We twist as mentioned above, and the proof of \cref{properness} goes through verbatim, to obtain a proper solution $\rho$ with $\# \rho(I_{\ell}) = 13$ for an auxiliary prime $\ell$. \cref{inertialorderdividesKg} again implies that $\rho$ does not arise from the $3$-torsion of an abelian threefold over $\Q$.
\end{proof}

\section{Acknowledgement}
The author was supported by Simons Foundation grant 550033. The author would like to thank Frank Calegari for helpful discussions, and Vishal Arul for valuable feedback that helped improve this paper.

\printbibliography

\end{document}